\documentclass[12pt,twoside,a4paper]{amsart}
\usepackage{amsmath,amssymb, mathtools}
\usepackage{mathtools}
\usepackage{cite}
\usepackage{graphicx}
\usepackage{txfonts}
\usepackage{mathrsfs}
\usepackage{empheq}
\newcommand{\va}{\varphi}
\newcommand{\p}{\psi}
\newcommand{\con}[3]{ (#1 \circ #2) (#3)}
\newcommand{\me}[2]{ {#1}_M ( #2 ) }

\allowdisplaybreaks[4]
\theoremstyle{definition}
\newtheorem{definition}{Definition}[subsection]
\newtheorem{remark}[definition]{Remark}

\theoremstyle{theorem}
\newtheorem{theorem}{Theorem}[subsection]

\newtheorem{fact}[theorem]{Fact}

\numberwithin{equation}{subsection}
\mathtoolsset{showonlyrefs=true}
\mathtoolsset{showonlyrefs=false}
\begin{document}

\title[Examples of Bessel Functions $\ldots$]{Examples of Bessel Functions of Order Zero Related to the Functional Equation of the Riemann Zeta Function}
\author{Hideto Iwata}
\address{ National Institute of Technology, Gunma College, Japan.}
\email{iwata@gunma-ct.ac.jp}
\maketitle
\vspace{0.5cm}

\textbf{Abstract.} 

In 2010, M.A.Chaudhry et al. introduced a class of functions $\Sigma$ satisfying an $L^1$-condition. Also, they proved that every ${\va} \in \Sigma$ leads to the functional equation for the Riemann zeta function. However, they also provided an example of a function outside $\Sigma$ for which the functional equation still holds.In this paper, we present several examples of Bessel functions of order zero that do not belong to $\Sigma$ but nevertheless satisfy the functional equation for the Riemann zeta function.In particular, we show that $J_0 (a\sqrt{x}), J_0 (ax)$, and $Y_0 (x)$ give rise to the functional equation despite not belonging to $\Sigma$.

\section{Introduction}

\subsection{Previous research}
The Riemann zeta function is defined by

\vspace{-0.3cm}

\begin{equation}
   \zeta(s) = \sum_{n=1}^\infty \frac{1}{n^s} = \prod_{p \text{: prime}} \left( 1 - \frac{1}{p^s} \right)^{-1} \quad (s=\sigma+it, \sigma>1).    \label{eq:1.1.1}
\end{equation}The function $\zeta(s)$ admits a meromorphic continuation to the whole complex plane and satisfies the following functional equation

\vspace{-0.3cm}

\begin{equation}
   \frac{\zeta(s)}{\zeta(1-s)}  = \pi^{s-1/2}\frac{\Gamma((1-s)/2)}{\Gamma(s/2)} = \chi(s),   \label{eq:1.1.2}
\end{equation}where

\vspace{-0.3cm}

\begin{align}
   \chi(s) 
   &= 2(2\pi)^{-(1-s)}\cos \left( \frac{\pi}{2}(1-s) \right)\Gamma(1-s)   \notag
   \\
   &= 2(2\pi)^{s-1}\sin \left( \frac{s}{2}\pi \right)\Gamma(1-s).       \label{eq:1.1.3}
   \end{align}The functional equation \eqref{eq:1.1.2} determines the locations of the zeros of $\zeta(s)$ and plays an important role in the study of prime numbers. By Stirling's formula, the order estimate of $\chi(s)$ is

\vspace{-0.3cm}

\begin{equation}
   \chi(1-s) \ll \left| \frac{t}{2\pi} \right|^{\sigma - \frac{1}{2}}   \label{eq:1.1.4}
\end{equation}as $|t| \to \infty$ for $0 < a \leq \sigma \leq b < 1/2$ and so $\chi(1-s)$ is analytic in $0 < \sigma < 1/2$. By \eqref{eq:1.1.4}, $\chi(1-s)$ converges to zero uniformly in the strip $0 < c_1 \leq c \leq c_2 < 1/2$ as $|t| \to \infty$. We define the function $\lambda(x)$ as the following complex integral :

\vspace{-0.3cm}

\begin{equation}
   \lambda(x) = \frac{1}{2\pi i}\int_{(c)} \chi(1-s)x^{-s}ds \quad (0 < c_1 \leq c \leq c_2 < 1/2).  \label{eq:1.1.5}
\end{equation}Here the path of integration in \eqref{eq:1.1.5} is the vertical line Re$(s) = c$ on the $s$-plane. By the estimate \eqref{eq:1.1.4}, the integral on the right hand side of \eqref{eq:1.1.5} is well defined for all $x > 0$. Moreover, a classical result in ~\cite{pk} yields $\lambda(x) = 2\cos (2\pi x)$. In ~\cite{ac}, M. A. Chaudhry et al. defined the convolution product (if it exists)

\vspace{-0.3cm}

\begin{equation}
   \con{ {\va} }{ {\p} }{ t } = \int_{0}^\infty {\va}(tx){\p}(x)dx.  \label{eq:1.1.6}
\end{equation}This product is noncommutative and satisfies the Mellin transform relation

\vspace{-0.3cm}

\begin{equation}
   M[ {\va} \circ {\p} ; s ] = \me{{\va}}{s}\me{ {\p} }{1-s},  \label{eq:1.1.7} 
\end{equation}where $\sigma_1 < \sigma < \sigma_2$ is a common strip of analyticity of both Mellin transforms

\vspace{-0.3cm}

\begin{gather*}
   \me{{\va}}{s} = M[ {\va}(x) ; s ] = \int_{0}^\infty {\va}(x) x^{s-1}dx,
   \\
   \me{{\p}}{1-s} = M[ {\p}(x) ; 1 - s ] = \int_{0}^\infty {\p}(x) x^{-s}dx.
\end{gather*}Also, M. A. Chaudhry et al. defined a class of functions $\Sigma$ as follows : For all positive $t$, a function ${\va}(t)$ is said to be a member of the class $\Sigma$ if $e^{bt}{\va}(t) \in L^1 (0,\infty)$ with $0 < b < 1$. For ${\va} \in \Sigma$, the convolution \eqref{eq:1.1.6} is well defined by Lemma 1 of ~\cite{ac}. Moreover, the Mellin transform $\me{{\va}}{s}$ exists by Lemma 2 of ~\cite{ac}. In ~\cite{ac}, every member of the class $\Sigma$ leads to the functional equation of the Riemann zeta function \eqref{eq:1.1.2} : 

\begin{theorem}[Theorem 1. in ~\cite{ac}]
For every ${\va} \in \Sigma$, we have the Riemann functional equation
\begin{equation}
   \me{\lambda}{s} = \frac{\zeta(1-s)}{\zeta(s)} = \frac{\me{{\p}}{s}}{\me{{\va}}{1-s}} \quad (0 <\sigma < 1/2),   \label{eq:1.1.8}
\end{equation}where ${\p} = \lambda \circ {\va}$.
\end{theorem}The condition ${\va}\in\Sigma$ is only a sufficient condition for the validity of \eqref{eq:1.1.8}, and it is not necessary. Such an example was given in ~\cite{ac} as follows: We consider the Bessel function of the first kind of order zero ${\va}(x) = J_0 (2\sqrt{x})$. Using \eqref{eq:3.2.1} with $a = 2$ and $\nu = 0$ below and substituting $u=\sqrt{x}$, we obtain $\me{{\va}}{s} = \Gamma(s)/\Gamma(1-s)$ for $0 < \text{Re}(s) < 3/4$ by \eqref{eq:3.2.1}. Next, we consider the function ${\p}(t)$ defined by the convolution product \eqref{eq:1.1.6}

\vspace{-0.5cm}

\begin{align}
   {\p}(t) 
   &= \con{\lambda}{{\va}}{t} = \int_{0}^\infty \lambda(tx){\va}(x)dx  \notag
   \\
   &= 2\int_{0}^\infty J_0 \left( 2\sqrt{x} \right)\cos(2\pi tx)dx = 2\left( \frac{1}{2\pi t}\sin\left( \frac{1}{2\pi t} \right) \right)   \label{eq:1.1.9}
\end{align}for all positive $t$. Here, we use the Fourier cosine transform \eqref{eq:3.1.1} with $a = 2$ and $\nu = 0$ below. Next, we calculate the Mellin transform for ${\p}(t)$. Substituting $u=t^{-1}$ and using \eqref{eq:3.2.2} with $a = 1/2\pi$ below, we have $\me{ {\p} }{s} = 2(2\pi)^{-s} \Gamma(1-s)\cos(\pi s/2)$ for $0 < \text{Re}(s) < 2$. Calculating the ratio ${\me{\va}{1-s}} / \me{{\p}}{s}$, we have the following relation 
\begin{equation}
   \frac{\me{{\va}}{1-s}}{\me{{\p}}{s}} = 2^{s-1} \pi^s \frac{1}{\Gamma(s)\cos\left( \frac{\pi s}{2} \right)}    \label{eq:1.1.10}
\end{equation}for $1/4 < \text{Re}(s) < 1$ Since the right hand side of \eqref{eq:1.1.10} is a meromorphic function, it is equal to $\zeta(s)/\zeta(1-s)$ by the functional equation for the Riemann zeta function. Therefore, the functional equation

\vspace{-0.3cm}

\begin{equation} 
   2^{s-1} \pi^s \frac{1}{\Gamma(s)\cos\left( \frac{\pi s}{2} \right)} = \frac{\zeta(s)}{\zeta(1-s)}    \label{eq:1.1.11}
\end{equation}holds for all complex numbers $s$ by analytic continuation.In this paper, we observe the same phenomenon for several other Bessel functions. In particular, all Bessel functions considered in this paper are of order zero. It is not clear at present why these examples lead to the same functional equation. Nevertheless, the results indicate that the phenomenon observed in Example 3 of ~\cite{ac} also occurs for several other Bessel functions of order zero.

\section{Main theorem}
\subsection{Main theorem}

We consider the following examples of Bessel functions of order zero.

\begin{theorem}
The following Bessel functions of order zero give rise to the functional equation of the Riemann zeta function through the ratio $\me{\va}{1-s} / \me{{\p}}{s}$:

\[  J_0 (a\sqrt{x}), J_0 (ax), Y_0 (x). \]Here, $Y_0 (x)$ denotes the Bessel function of the second kind of order zero, and $a$ is a positive real constant. 
\end{theorem}

The example $J_0(2\sqrt{x})$ was considered in~\cite{ac}. Theorem 2.1.1 shows that the same phenomenon persists for the more general family $J_0(a\sqrt{x})$ for all positive real numbers $a$. For the Bessel function, some of these functions belong to $\Sigma$, while others do not. For example, the modified Bessel function of the third kind of order zero $K_0$ belongs to $\Sigma$ since the following order estimates hold (see~\cite{L}): 

\vspace{-0.3cm}

\begin{gather}
   K_0 (x) \ll x^{-\frac{1}{2}}e^{-x} \quad (x \to \infty),   \label{eq:2.1.1}
   \\
   K_0 (x) \sim \log \left( \frac{2}{x} \right) \quad (x \to 0).   \label{eq:2.1.2}
\end{gather}Hence, $K_0$ leads to the functional equation by Theorem 1.1.1. On the other hand, $J_0$ does not belong to $\Sigma$ since the following order estimate holds (see~\cite{L}): 

\vspace{-0.3cm}

\begin{equation}
   J_0 (x) \ll x^{-\frac{1}{2}}  \quad (x \to \infty),   \label{eq:2.1.3}
\end{equation}

\section{Preliminary}
\subsection{The Fourier cosine transform}
In this subsection, we collect the classical results of the Fourier cosine transform for the $J,Y$-Bessel functions when we prove Theorem 2.1.1. Throughout this paper, $a$ denotes a positive real constant.

\begin{fact}[~\cite{b}]
\begin{gather}
     \int_{0}^\infty J_0 \left( a\sqrt{x} \right)\cos (yx)dx
     = y^{-1}\sin\left( \frac{a^2}{4}y^{-1} \right)  \quad (y>0)   \label{eq:3.1.1}
     \\
     \int_{0}^\infty J_0 (ax)\cos (yx)dx 
  = \begin{cases}
      (a^2 - y^2)^{-1/2} \quad (0 < y < a)
      \\
      \infty \quad (y = a)
      \\
      0 \quad (a < y < \infty)    \label{eq:3.1.2}
     \end{cases}
     \\     
     \int_{0}^\infty Y_0 (ax)\cos (yx)dx 
  = \begin{cases}
      0 \quad (0 < y < a)
      \\
      -(y^2 - a^2)^{-1/2} \quad (a < y < \infty)   \label{eq:3.1.3}
     \end{cases}
\end{gather}
\end{fact}

\subsection{The Mellin transform}
Next, we collect the classical results of the Mellin transform when we prove Theorem 2.1.1.  

\begin{fact}[~\cite{b}]
\begin{gather}
    \int_{0}^\infty J_\nu (ax)x^{s-1}dx = 2^{s-1}a^{-s} \frac{\Gamma\left( \frac{s}{2} + \frac{\nu}{2} \right)}{\Gamma\left( \frac{\nu}{2} - \frac{s}{2} + 1 \right)}  \quad \left(a > 0, -\text{Re}(\nu) < \text{Re}(s) < \frac{3}{2} \right),   \label{eq:3.2.1}
    \\
    \int_{0}^\infty \sin(ax) x^{-s}dx = a^{s-1}\Gamma(1-s)\sin\left( \frac{\pi}{2}(1 - s) \right) \quad (a>0, 0 < \text{Re}(s) < 2),  \label{eq:3.2.2}
    \\
    \int_{0}^\infty Y_0 (ax)x^{s-1}dx = -2^{s-1} a^{-s} \pi^{-1} \left\{ \Gamma\left( \frac{s}{2} \right) \right\}^2 \cos \left( \frac{\pi}{2}s \right) \quad \left(a>0, 0 <  \text{Re}(s) < \frac{3}{2} \right),  \label{eq:3.2.3}    
\end{gather}
\end{fact}

\section{Proof of Theorem 2.1.1}

\subsection{The case ${\va}(x) = J_0 (a\sqrt{x})$}
First, we prove Theorem 2.1.1 for the case ${\va}(x) = J_0 (a\sqrt{x})$ and calculate the following Mellin transform

\vspace{-0.3cm}

\begin{equation} 
   \me{{\va}}{s} = \int_{0}^\infty J_0 (a\sqrt{x})x^{s-1}dx.  \label{eq:4.1.1} 
\end{equation}With the substitution $t = a\sqrt{x}$ and using \eqref{eq:3.2.1} with $\nu=0$ and $a=1$, we have 

\vspace{-0.3cm}

\begin{equation}
   \me{{\va}}{s} = a^{-2s}2^{2s}\frac{\Gamma(s)}{\Gamma(1-s)}   \label{eq:4.1.2}
\end{equation}for $0 < \text{Re}(s) < 3/4$. Next we consider the function ${\p}(t)$ as the same in \eqref{eq:1.1.9} ;

\vspace{-0.3cm}

\begin{equation}
   {\p}(t) = 2\int_{0}^\infty J_0 (a\sqrt{x})\cos(2\pi tx)dx   \label{eq:4.1.3} 
\end{equation}for all positive $t$. Using \eqref{eq:3.1.1} with $y = 2\pi t$, we obtain ${\p}(t) = (\pi t)^{-1}\sin\left( a^2 / 8\pi t \right)$. Next, we calculate the Mellin transform for ${\p}(t)$, that is

\vspace{-0.3cm}

\begin{equation}
   \me{\p}{s} = \frac{1}{\pi}\int_{0}^\infty \frac{1}{x}\sin\left( \frac{a^2}{8\pi x} \right) x^{s-1}dx.   \label{eq:4.1.4} 
\end{equation}With the substitution $u = x^{-1}$ in \eqref{eq:4.1.4} and using \eqref{eq:3.2.2}, we have

\vspace{-0.3cm}

\begin{equation}
   \me{\p}{s} = 2^{3-3s} a^{2s-2}\pi^{-s} \Gamma(1-s)\cos \left( \frac{s}{2}\pi \right)   \label{eq:4.1.5}
\end{equation}for $0 < \text{Re}(s) < 2$. Calculating the ratio ${\me{\va}{1-s}} / \me{{\p}}{s}$, we have the relation
 
\vspace{-0.3cm}

\begin{equation}
   \frac{\me{{\va}}{1-s}}{\me{{\p}}{s}} = 2^{s-1} \pi^s \frac{1}{\Gamma(s)\cos \left( \frac{s}{2}\pi \right)}    \label{eq:4.1.6}
\end{equation}for $0 < \text{Re}(s) < 1$. Also, we use a fundamental property for $\Gamma$-function, we can rewrite \eqref{eq:4.1.6} as

\begin{equation}  
   \frac{\me{{\va}}{1-s}}{\me{{\p}}{s}} = 2(2\pi)^{s-1} \Gamma(1-s)\sin\left( \frac{s}{2}\pi \right).   \label{eq:4.1.7}
\end{equation}By the functional equation \eqref{eq:1.1.2} and \eqref{eq:1.1.3}, the right hand side equals $\zeta(s) / \zeta(1-s)$, and we have the relation  

\vspace{-0.3cm}

\begin{equation}
   2(2\pi)^{s-1} \Gamma(1-s)\sin\left( \frac{s}{2}\pi \right) = \frac{\zeta(s)}{\zeta(1-s)}.   \label{eq:4.1.8}
\end{equation}Hence, \eqref{eq:4.1.8} holds for all complex numbers $s$ by analytic continuation. Therefore, Theorem 2.1.1 is proved in this case.  

\subsection{The case ${\va}(x) = J_0 (ax)$}
Secondly, we prove Theorem 2.1.1 for the case ${\va}(x) = J_0 (ax)$ and $\me{{\va}}{s} = a^{-s}2^{s-1}\Gamma(s/2) / \Gamma(1-s/2)$ for $0 < \text{Re}(s) < 3/2$ by \eqref{eq:3.2.1} with $\nu = 0$. Next we consider the function ${\p}(t)$ as the same in \eqref{eq:1.1.9} ;

\vspace{-0.3cm}

\begin{equation}
   {\p}(t) = 2\int_{0}^\infty J_0 (ax)\cos(2\pi tx)dx   \label{eq:4.2.1} 
\end{equation}for all positive $t$. Using \eqref{eq:3.1.2} with $y = 2\pi t$, we obtain ${\p}(t) = 2(a^2 -4\pi^2 t^2)^{-1/2}$ for $0 < t < a/2\pi$. Next, we calculate the Mellin transform for ${\p}(t)$, that is

\vspace{-0.3cm}

\begin{equation}
   \me{\p}{s} = 2\int_{0}^{\frac{a}{2\pi}} (a^2 - 4\pi^2 x^2)^{-\frac{1}{2}}x^{s-1}dx.   \label{eq:4.2.2} 
\end{equation}With the substitution $u = 2\pi x/a$ in \eqref{eq:4.2.2} and using the following classical Beta integral formula

\[ B(x,y) = \int_{0}^1 u^{x-1} (1-u)^{y-1}du = \frac{\Gamma(x)\Gamma(y)}{\Gamma(x+y)} \] for $\text{Re}(x), \text{Re}(y)>0$, we obtain

\vspace{-0.3cm}

\begin{align*}
   \me{\p}{s}
   &= \frac{2}{a}\left( \frac{a}{2\pi} \right)^s \int_{0}^1 u^{s-1}(1-u^2)^{-\frac{1}{2}}du
   \\
   &=  \frac{2}{a}\left( \frac{a}{2\pi} \right)^s \frac{1}{2}B\left( \frac{s}{2}, \frac{1}{2} \right)
    \\
    &= a^{s-1} \pi^{\frac{1}{2} - s} 2^{-s} \frac{\Gamma(s/2)}{\Gamma((s+1)/2)}
\end{align*}for $\text{Re}(s) > 0$. Calculating the ratio ${\me{\va}{1-s}} / \me{{\p}}{s}$, we have

\vspace{-0.3cm}

\begin{equation}
   \frac{\me{{\va}}{1-s}}{\me{{\p}}{s}} = \pi^{s - \frac{1}{2}} \frac{\Gamma((1-s)/2)}{\Gamma(s/2)}  \label{eq:4.2.3}
\end{equation}for $0< \text{Re}(s) < 1$. By the same argument as above, the right hand side of \eqref{eq:4.2.3} is equal to $\zeta(s)/\zeta(1-s)$. Hence, the functional equation holds for all complex numbers $s$ by analytic continuation. Therefore, Theorem 2.1.1 is proved in this case.  

\subsection{The case ${\va}(x) = Y_0 (x)$}
Finally, we prove Theorem 2.1.1 for the case ${\va}(x) = Y_0 (x)$ and $\me{{\va}}{s} = -2^{s-1}\{\Gamma(s/2)\}^2 \cos(s\pi/2)$ for $0 < \text{Re}(s) < 3/2$ by \eqref{eq:3.2.3} with $a=1$. Next we consider the function ${\p}(t)$ as the same in \eqref{eq:1.1.9} ;

\vspace{-0.3cm}

\begin{equation}
   {\p}(t) = 2\int_{0}^\infty Y_0 (x)\cos(2\pi tx)dx   \label{eq:4.3.1} 
\end{equation}for all positive $t$. Using \eqref{eq:3.1.3} with $y=2\pi t$, we obtain ${\p}(t) = -2(4\pi^2 t^2 - 1)^{-1/2}$ for $t > 1/2\pi$. Next, we calculate the Mellin transform for ${\p}(t)$, that is

\vspace{-0.3cm}

\begin{equation}
   \me{\p}{s} = -2\int_{\frac{1}{2\pi}}^\infty ( 4\pi^2 x^2 - 1 )^{-\frac{1}{2}}x^{s-1}dx.   \label{eq:4.3.2} 
\end{equation}With the substitution $u = 2\pi x$ in \eqref{eq:4.3.2} and using the following classical Beta integral formula

\[ \int_{0}^\frac{\pi}{2} \sin^{x-1}\theta \cos^{y-1}\theta d\theta = \frac{1}{2}B\left( \frac{x}{2}, \frac{y}{2} \right)  \quad \left( \text{Re}(x), \text{Re}(y) > 0 \right), \] we have

\vspace{-0.3cm}

\begin{equation}
   \me{\p}{s} = -2^{-s}\pi^{-s+\frac{1}{2}}\frac{\Gamma( (1-s) / 2 )}{\Gamma( 1 - s/2 )}   \label{eq:4.3.3}
\end{equation}for $\text{Re}(s) < 1$. Calculating the ratio ${\me{\va}{1-s}} / \me{{\p}}{s}$ and using a fundamental property for $\Gamma$-function, we have

\vspace{-0.3cm}

\begin{align}
   \frac{\me{{\va}}{1-s}}{\me{{\p}}{s}}   \notag
   &= \pi^{s - \frac{3}{2}} \Gamma\left( \frac{1-s}{2} \right)\Gamma\left(1 - \frac{s}{2} \right)\sin\left( \frac{s}{2}\pi \right)    \notag
   \\
   &= 2(2\pi)^{s-1} \Gamma(1-s)\sin\left( \frac{s}{2}\pi \right)   \label{eq:4.3.4}
\end{align}for $-1/2 < \text{Re}(s) < 1$. By the functional equation \eqref{eq:1.1.2} and \eqref{eq:1.1.3}, the expression in \eqref{eq:4.3.4} is equal to $\zeta(s)/\zeta(1-s)$. Hence, the functional equation holds for all complex numbers $s$ by analytic continuation. Therefore, Theorem 2.1.1 is proved in this case. Combining the three cases above, we obtain Theorem 2.1.1.   $\square$

\begin{remark}
The modified Bessel function of the first kind of order zero $I_0$ also does not belong to $\Sigma$, and the above method based on the Fourier and Mellin transforms cannot be applied directly. The reason is that the corresponding integral transforms do not converge, since the following asymptotic formula holds for $I_0$ (see ~\cite{L}) :

\vspace{-0.3cm}

\begin{equation}
   I_0 (x) \sim \frac{e^x}{\sqrt{2\pi x}} \quad (x \to \infty).   \label{eq:4.3.5}
\end{equation}

\end{remark}

\end{document}